\documentclass[12pt,thmsa]{article}%
\usepackage{sw20bams}
\usepackage{amsmath}
\usepackage{amsfonts}
\usepackage{amssymb}
\usepackage{graphicx}%
\providecommand{\U}[1]{\protect\rule{.1in}{.1in}}
\begin{document}

\author{Steven R. Finch}
\title{Correlation between Angle and Side}
\date{March 21, 2022}
\maketitle

\begin{abstract}
Let $\alpha$ be an arbitrary angle in a random spherical triangle $\Delta$ and
$a$ be the side opposite $\alpha$. (The sphere has radius $1$; vertices of
$\Delta$ are independent and uniform.) If some other side is constrained to be
$\pi/2$, then $\operatorname*{E}(\alpha\,a)=3.05...$. If instead some other
angle is fixed at $\pi/2$, then $\operatorname*{E}(\alpha\,a)=2.87...$. In our
study of the latter scenario, both Ap\'ery's constant and Catalan's constant
emerge. We also review Miles' 1971 proof that $\operatorname*{E}%
(\alpha\,a)=\pi^{2}/2-2$ when no constraints are in place.

\end{abstract}

\footnotetext{Copyright \copyright \ 2010, 2017, 2022 by Steven R. Finch. All
rights reserved.}For any planar triangle, long sides are opposite large angles
and short sides are opposite small angles. Quantifying this observation for
random triangles with either

\begin{itemize}
\item independent Gaussian vertices in the plane, or

\item independent uniform vertices in a compact convex subset of the plane
\end{itemize}

\noindent seems analytically intractable. We turn attention therefore to
random spherical triangles with independent uniform vertices on the unit sphere.

A\ spherical triangle $\Delta$ is a region enclosed by three great circles on
the sphere; a great circle is a circle whose center is at the origin. The
sides of $\Delta$ are arcs of great circles and have length $a$, $b$, $c$.
Each of these is $\leq\pi$. The angle $\alpha$ opposite side $a$ is the
dihedral angle between the two planes passing through the origin and
determined by arcs $b$, $c$. The angles $\beta$, $\gamma$ opposite sides $b $,
$c$ are similarly defined. Each of these is $\leq\pi$ too.

Given a random spherical triangle, the univariate density for $a$ is
\[%
\begin{array}
[c]{ccc}%
\dfrac12\sin(a), &  & 0<a<\pi
\end{array}
\]
and
\[%
\begin{array}
[c]{ccc}%
\operatorname*{E}(a)=\dfrac\pi2, &  & \operatorname*{E}(a^{2})=\dfrac{\pi^{2}%
}2-2.
\end{array}
\]
Further, $\alpha$ is uniformly distributed on $[0,\pi]$ and
\[%
\begin{array}
[c]{ccc}%
\operatorname*{E}(\alpha)=\dfrac\pi2, &  & \operatorname*{E}(\alpha
^{2})=\dfrac{\pi^{2}}3.
\end{array}
\]
It can be shown that $\alpha$, $b$, $c$ are independent random variables;
hence $\operatorname*{E}(\alpha\,b)=\pi^{2}/4=\operatorname*{E}(\alpha\,c)$.
In contrast, the density for $(a,\beta,\gamma)$ is \cite{Ms}
\[
\frac1{4\pi}\frac{\sin(\beta)\sin(\gamma)\sin(a)^{3}}{\left(  1-(\cos
(\beta)\cos(\gamma)-\sin(\beta)\sin(\gamma)\cos(a))^{2}\right)  ^{3/2}}.
\]
As special cases, the conditional density for $(\beta,\gamma)$ given that
$a=\pi/2$ is
\[
\frac1{2\pi}\frac{\sin(\beta)\sin(\gamma)}{\left(  1-\cos(\beta)^{2}%
\cos(\gamma)^{2}\right)  ^{3/2}};
\]
the conditional density for $(a,\gamma)$ given that $\beta=\pi/2$ is
\[
\frac14\frac{\sin(\gamma)\sin(a)^{3}}{\left(  1-\sin(\gamma)^{2}\cos
(a)^{2}\right)  ^{3/2}};
\]
and the unconditional density for $(\beta,\gamma)$ is
\[
\frac1{2\pi}\frac1{\sin(\beta)^{2}\sin(\gamma)^{2}}\cdot\left\{
\begin{array}
[c]{lll}%
-\cos(\gamma)\sin(\gamma)+\gamma &  & \text{if }\beta-\gamma>0\text{ and
}\beta+\gamma<\pi\text{,}\\
\pi+\cos(\gamma)\sin(\gamma)-\gamma &  & \text{if }\beta-\gamma<0\text{ and
}\beta+\gamma>\pi\text{,}\\
-\cos(\beta)\sin(\beta)+\beta &  & \text{if }\beta-\gamma<0\text{ and }%
\beta+\gamma<\pi\text{,}\\
\pi+\cos(\beta)\sin(\beta)-\beta &  & \text{if }\beta-\gamma>0\text{ and
}\beta+\gamma>\pi.
\end{array}
\right.
\]
These facts will be needed later.

\section{Univariate Densities}

Sides $a$, $b$, $c$ are pairwise independent; thus the conditional density for
$b$ given $c=\pi/2$ remains unchanged (the sine density on $[0,\pi]$). Angles
$\alpha$, $\beta$, $\gamma$ are uncorrelated but pairwise \textit{dependent}.
Therefore the case of two angles, plus two other scenarios involving opposite
side and angle, yield interesting new results.

\subsection{Angle $\beta,$ for Fixed Angle $\gamma$}

The conditional density for $\beta$ given that $\gamma=\pi/2$ is
\begin{align*}
&  \frac12\frac1{\sin(\beta)^{2}}\cdot\left\{
\begin{array}
[c]{lll}%
-\cos(\beta)\sin(\beta)+\beta &  & \text{if }0<\beta<\pi/2\text{,}\\
\pi+\cos(\beta)\sin(\beta)-\beta &  & \text{if }\pi/2<\beta<\pi
\end{array}
\right. \\
&  =\frac12\cdot\left\{
\begin{array}
[c]{lll}%
-\cot(\beta)+\beta\csc(\beta)^{2} &  & \text{if }0<\beta<\pi/2\text{,}\\
\cot(\beta)+(\pi-\beta)\csc(\beta)^{2} &  & \text{if }\pi/2<\beta<\pi.
\end{array}
\right.
\end{align*}
It follows that
\[%
\begin{array}
[c]{ccc}%
\operatorname*{E}\left(  \beta\left|  \gamma=\dfrac\pi2\right.  \right)
=\dfrac\pi2, &  & \operatorname*{E}\left(  \beta^{2}\left|  \gamma=\dfrac
\pi2\right.  \right)  =\dfrac{\pi^{2}}2-\dfrac74\zeta(3)
\end{array}
\]
where
\[
\zeta(3)=%
{\displaystyle\sum\limits_{k=1}^{\infty}}
\frac1{k^{3}}
\]
is Ap\'ery's constant \cite{F1}.

\subsection{Side $c$, for Fixed Angle $\gamma$}

By the Law of Cosines for Sides:
\[
\cos(c)=\cos(a)\cos(b)+\sin(a)\sin(b)\cos(\gamma)
\]
we obtain
\[
\cos(c)=\cos(a)\cos(b)
\]
if $\gamma=\pi/2$. Let $u=\cos(a)$, $v=\cos(b)$, $w=u\,v$, $z=\arccos(w) $.
Then $u$, $v$ are independent uniform on $[-1,1]$, that is, with density
\[
f(u,v)=\left\{
\begin{array}
[c]{lll}%
1/4 &  & \text{if }-1\leq u\leq1\text{ and }-1\leq v\leq1,\\
0 &  & \text{otherwise.}%
\end{array}
\right.
\]
By \cite{Pa, GLD}, the density of $w$ is
\[
g(w)=%
{\displaystyle\int\limits_{-\infty}^{\infty}}
f\left(  t,\frac wt\right)  \frac1{|t|}dt=\frac14%
{\displaystyle\int\limits_{-1}^{1}}
\varepsilon(w,t)\frac1{|t|}dt
\]
where $\varepsilon(w,t)=1$ if $-1<w/t<1$, $\varepsilon(w,t)=0$ otherwise. We
obtain
\[
g(w)=-\frac12\ln|w|.
\]
Since $0\leq z\leq\pi$ and
\[
\left|  \frac{dz}{dw}\right|  =\frac1{\sqrt{1-w^{2}}}=\frac1{\sin(z)},
\]
the density of $z$ is
\[
h(z)=\frac{g(\cos(z))}{\frac1{\sin(z)}}=-\frac12\sin(z)\ln|\cos(z)|.
\]
It follows that the conditional density for side $c$, given $\gamma=\pi/2$,
has a singularity at $c=\pi/2$ and
\[%
\begin{array}
[c]{ccc}%
\operatorname*{E}\left(  c\left|  \gamma=\dfrac\pi2\right.  \right)
=\dfrac\pi2, &  & \operatorname*{E}\left(  c^{2}\left|  \gamma=\dfrac
\pi2\right.  \right)  =-6+\dfrac{\pi^{2}}2+4G
\end{array}
\]
where
\[
G=%
{\displaystyle\sum\limits_{k=0}^{\infty}}
\frac{(-1)^{k}}{(2k+1)^{2}}
\]
is Catalan's constant \cite{F2}.

\subsection{Angle $\gamma,$ for Fixed Side $c$}

By the Law of Cosines for Angles:
\[
-\cos(\gamma)=\cos(\alpha)\cos(\beta)-\sin(\alpha)\sin(\beta)\cos(c)
\]
we obtain
\[
\cos(\gamma)=-\cos(\alpha)\cos(\beta)
\]
if $c=\pi/2$. Let $u=\cos(\alpha)$, $v=\cos(\beta)$, $w=-u\,v$, $z=\arccos
(w)$. The Jacobian determinant of $(\alpha,\beta)\mapsto(u,v)$ is
\[
\left|
\begin{array}
[c]{cc}%
-\sin(\alpha) & 0\\
0 & -\sin(\beta)
\end{array}
\right|  =\sin(\alpha)\sin(\beta)=\sqrt{1-u^{2}}\sqrt{1-v^{2}}
\]
because $0\leq\alpha\leq\pi$, $0\leq\beta\leq\pi$. Thus $u$, $v$ have density
\[
\frac1{2\pi}\frac{\sqrt{1-u^{2}}\sqrt{1-v^{2}}}{\left(  1-u^{2}v^{2}\right)
^{3/2}}\frac1{\sqrt{1-u^{2}}\sqrt{1-v^{2}}}=\frac1{2\pi}\frac1{\left(
1-u^{2}v^{2}\right)  ^{3/2}}.
\]
By \cite{Pa, GLD}, the density of $w$ is
\[
g(w)=%
{\displaystyle\int\limits_{-\infty}^{\infty}}
f\left(  t,\frac wt\right)  \frac1{|t|}dt=\frac1{2\pi}\frac1{\left(
1-w^{2}\right)  ^{3/2}}%
{\displaystyle\int\limits_{-1}^{1}}
\varepsilon(w,t)\frac1{|t|}dt
\]
where $\varepsilon(w,t)=1$ if $-1<w/t<1$, $\varepsilon(w,t)=0$ otherwise. We
obtain
\[
g(w)=-\frac1\pi\frac{\ln|w|}{\left(  1-w^{2}\right)  ^{3/2}}
\]
and hence the density of $z$ is
\[
h(z)=\frac{g(\cos(z))}{\frac1{\sin(z)}}=-\frac1\pi\frac{\ln|\cos(z)|}%
{\sin(z)^{2}}.
\]
It follows that the conditional density for angle $\gamma$, given $c=\pi/2$,
has a singularity at $\gamma=\pi/2$ and
\[%
\begin{array}
[c]{ccc}%
\operatorname*{E}\left(  \gamma\left|  c=\dfrac\pi2\right.  \right)
=\dfrac\pi2, &  & \operatorname*{E}\left(  \gamma^{2}\left|  c=\dfrac
\pi2\right.  \right)  =\dfrac{\pi^{2}}4+\ln(2)^{2}.
\end{array}
\]
This completes our quick survey of univariate densities, for a fixed side or angle.

\section{Bivariate Moments}

We evaluate $\operatorname*{E}(\alpha\,a\,|\,b=\pi/2)$ and $\operatorname*{E}%
(\alpha\,a\,|\,\beta=\pi/2)$ here, giving precise numerics for the former and
exact symbolics for the latter.

\subsection{(Angle $\alpha$, Side $a)$, for Fixed Side $b$}

The Law of Cosines for Sides:
\[
\cos(a)=\cos(b)\cos(c)+\sin(b)\sin(c)\cos(\alpha)
\]
can be expressed as
\[
w=u\,v+\sqrt{1-u^{2}}\sqrt{1-v^{2}}\cos(\theta)
\]
where $u=\cos(b)$, $v=\cos(c)$, $w=\cos(a)$, $\theta=\alpha$. Then $u$, $v$,
$\theta$ are independent; $u$, $v$, $w$ are uniform on $[-1,1]$ in the
unconditional case and $\theta$ is uniform on $[0,\pi]$. Fix $0\leq b\leq
\pi/2$ for simplicity, then $0\leq u\leq1$. Solving for $v$ in terms of $u$,
$w$, $\theta$ we obtain two solutions
\[
\varphi(u,w,\theta)=\frac{u\,w+|\cos(\theta)|\sqrt{\left(  1-u^{2}\right)
\left(  u^{2}-w^{2}+(1-u^{2})\cos(\theta)^{2}\right)  }}{u^{2}+(1-u^{2}%
)\cos(\theta)^{2}},
\]
\[
\psi(u,w,\theta)=\frac{u\,w-|\cos(\theta)|\sqrt{\left(  1-u^{2}\right)
\left(  u^{2}-w^{2}+(1-u^{2})\cos(\theta)^{2}\right)  }}{u^{2}+(1-u^{2}%
)\cos(\theta)^{2}}
\]
assuming
\[
u^{2}-w^{2}+(1-u^{2})\cos(\theta)^{2}>0
\]
and, further,
\[
\left(  w>-u\;\text{and\ }\theta<\pi/2\right)  \;\text{or\ }\left(
w<-u\;\text{and\ }\theta>\pi/2\right)
\]
for $\varphi$ and
\[
\left(  w>u\;\text{and\ }\theta<\pi/2\right)  \;\text{or\ }\left(
w<u\;\text{and\ }\theta>\pi/2\right)
\]
for $\psi$. Observe that the domains for $\varphi$, $\psi$ overlap when
\[
\left(  w>u\;\text{and\ }\theta<\pi/2\right)  \;\text{or\ }\left(
w<-u\;\text{and\ }\theta>\pi/2\right)  ,
\]
that is, the transformation is one-to-one for $(w,\theta)\in[-u,u]\times
[0,\pi]$ and two-to-one otherwise. Also, the Jacobian determinant of
$(v,\theta)\mapsto(w,\theta)$ is
\[
\delta(u,v,\theta)=u-\frac{\sqrt{1-u^{2}}v\cos(\theta)}{\sqrt{1-v^{2}}}.
\]
Let
\[
\xi(u,\theta)=\sqrt{u^{2}+(1-u^{2})\cos(\theta)^{2}}
\]
for convenience, then $\operatorname*{E}(\alpha\,a\,|\,b)$ is equal to
\cite{Pa}
\begin{align*}
&  \ \ \ \frac1{2\pi}%
{\displaystyle\int\limits_{0}^{\pi/2}}
\,%
{\displaystyle\int\limits_{-u}^{u}}
\frac{\theta\arccos(w)}{\left|  \delta(u,\varphi(u,w,\theta),\theta)\right|
}dw\,d\theta+\frac1{2\pi}%
{\displaystyle\int\limits_{\pi/2}^{\pi}}
\,%
{\displaystyle\int\limits_{-u}^{u}}
\frac{\theta\arccos(w)}{\left|  \delta(u,\psi(u,w,\theta),\theta)\right|
}dw\,d\theta\\
&  \ \ +\frac1{2\pi}%
{\displaystyle\int\limits_{0}^{\pi/2}}
\,%
{\displaystyle\int\limits_{u}^{\xi(u,\theta)}}
\left(  \frac1{\left|  \delta(u,\varphi(u,w,\theta),\theta)\right|  }%
+\frac1{\left|  \delta(u,\psi(u,w,\theta),\theta)\right|  }\right)
\theta\arccos(w)dw\,d\theta\\
&  \ \ +\frac1{2\pi}%
{\displaystyle\int\limits_{\pi/2}^{\pi}}
\,%
{\displaystyle\int\limits_{-\xi(u,\theta)}^{-u}}
\left(  \frac1{\left|  \delta(u,\varphi(u,w,\theta),\theta)\right|  }%
+\frac1{\left|  \delta(u,\psi(u,w,\theta),\theta)\right|  }\right)
\theta\arccos(w)dw\,d\theta.
\end{align*}
In the event $b=\pi/2$, we have $u=0$,
\[
\varphi(0,w,\theta)=\frac{\sqrt{-w^{2}+\cos(\theta)^{2}}}{|\cos(\theta
)|}=-\psi(0,w,\theta),
\]
\[%
\begin{array}
[c]{ccc}%
\delta(0,v,\theta)=-\dfrac{v\cos(\theta)}{\sqrt{1-v^{2}}}, &  & \xi
(0,\theta)=|\cos(\theta)|,
\end{array}
\]
\begin{align*}
\frac1{\left|  \delta(0,\varphi(0,w,\theta),\theta)\right|  }+\frac1{\left|
\delta(0,\psi(0,w,\theta),\theta)\right|  }  &  =\frac{\sqrt{1-\varphi^{2}}%
}{\varphi\,|\cos(\theta)|}+\frac{\sqrt{1-\psi^{2}}}{(-\psi)|\cos(\theta)|}\\
\  &  =\frac{2\sqrt{1-\varphi^{2}}}{\varphi\,|\cos(\theta)|}%
\end{align*}
which becomes
\[
\frac{2\sqrt{1-\frac{-w^{2}+\cos(\theta)^{2}}{\cos(\theta)^{2}}}}{\sqrt
{-w^{2}+\cos(\theta)^{2}}}=\frac{2|w|}{|\cos(\theta)|\sqrt{-w^{2}+\cos
(\theta)^{2}}}
\]
and therefore $\operatorname*{E}(\alpha\,a\,|\,b=\pi/2)$ is equal to
\[
\ \frac1\pi%
{\displaystyle\int\limits_{0}^{\pi/2}}
\,%
{\displaystyle\int\limits_{0}^{\cos(\theta)}}
\frac{\theta\,w\arccos(w)}{\cos(\theta)\sqrt{-w^{2}+\cos(\theta)^{2}}%
}dw\,d\theta+\frac1\pi%
{\displaystyle\int\limits_{\pi/2}^{\pi}}
\,%
{\displaystyle\int\limits_{\cos(\theta)}^{0}}
\frac{\theta\,w\arccos(w)}{\cos(\theta)\sqrt{-w^{2}+\cos(\theta)^{2}}%
}dw\,d\theta.
\]
This can be reduced to a single integral:
\[
\frac14%
{\displaystyle\int\limits_{0}^{\pi}}
\,\left[  2-\,_{2}F_{1}\left(  \tfrac12,\tfrac12,2,\cos(\theta)^{2}\right)
\cos(\theta)\right]  \theta\,d\theta=3.0538319164380270202505577...
\]
involving the following Gauss hypergeometric function:
\begin{align*}
_{2}F_{1}\left(  \tfrac12,\tfrac12,2,x\right)   &  =\frac1\pi%
{\displaystyle\sum\limits_{n=0}^{\infty}}
\frac{\Gamma(n+1/2)^{2}}{\Gamma(n+2)}\frac{x^{n}}{n!}\\
&  =\frac4\pi\left[  \frac1x%
{\displaystyle\int\limits_{0}^{\pi/2}}
\sqrt{1-x\sin(t)^{2}}\,dt+\left(  1-\frac1x\right)
{\displaystyle\int\limits_{0}^{\pi/2}}
\dfrac1{\sqrt{1-x\sin(t)^{2}}}\,dt\right]  .
\end{align*}
Despite a connection to complete elliptic integrals \cite{FJ}, this
unfortunately seems to be as far as we can go.

\subsection{(Angle $\alpha$, Side $a)$, for Fixed Angle $\beta$}

The Law of Cosines for Angles:
\[
-\cos(\alpha)=\cos(\beta)\cos(\gamma)-\sin(\beta)\sin(\gamma)\cos(a)
\]
can be expressed as
\[
w=-u\,v+\sqrt{1-u^{2}}\sqrt{1-v^{2}}\cos(\theta)
\]
where $u=\cos(\beta)$, $v=\cos(\gamma)$, $w=\cos(\alpha)$, $\theta=a$. Fix
$0\leq\beta\leq\pi/2$ for simplicity, then $0\leq u\leq1$. Solving for $v$ in
terms of $u$, $w$, $\theta$ we obtain two solutions $\varphi(-u,w,\theta)$,
$\psi(-u,w,\theta)$ as before. Also, the Jacobian determinant of
$(v,\theta)\mapsto(w,\theta)$ is
\[
\delta(u,v,\theta)=-u-\frac{\sqrt{1-u^{2}}v\cos(\theta)}{\sqrt{1-v^{2}}}.
\]
In the event $\beta=\pi/2$, we have $u=0$ and an identical formula for
$\left|  \delta(0,\varphi,\theta)\right|  ^{-1}+\left|  \delta(0,\psi
,\theta)\right|  ^{-1}$ follows. The distinction with earlier calculations
arises from the density
\[
\frac14\frac{\sqrt{1-v^{2}}\sin(\theta)^{3}}{\left(  1-(1-v^{2})\cos
(\theta)^{2}\right)  ^{3/2}}\frac1{\sqrt{1-v^{2}}}=\frac14\frac{\sin
(\theta)^{3}}{\left(  1-(1-v^{2})\cos(\theta)^{2}\right)  ^{3/2}}
\]
for $(v,\theta)$. Substituting $\varphi$ in place of $v$, we obtain
\[
\frac14\frac{\sin(\theta)^{3}}{\left(  1-\left[  1-\frac{-w^{2}+\cos
(\theta)^{2}}{\cos(\theta)^{2}}\right]  \cos(\theta)^{2}\right)  ^{3/2}}%
=\frac14\frac{\sin(\theta)^{3}}{\left(  1-w^{2}\right)  ^{3/2}}
\]
and therefore $\operatorname*{E}(\alpha\,a\,|\,\beta=\pi/2)$ is equal to
\begin{align*}
&  \ \ \ \ \ \frac12%
{\displaystyle\int\limits_{0}^{\pi/2}}
\,%
{\displaystyle\int\limits_{0}^{\cos(\theta)}}
\frac{\theta\,w\arccos(w)}{\cos(\theta)\sqrt{-w^{2}+\cos(\theta)^{2}}}%
\frac{\sin(\theta)^{3}}{\left(  1-w^{2}\right)  ^{3/2}}dw\,d\theta\\
&  \ \ \ \ +\frac12%
{\displaystyle\int\limits_{\pi/2}^{\pi}}
\,%
{\displaystyle\int\limits_{\cos(\theta)}^{0}}
\frac{\theta\,w\arccos(w)}{\cos(\theta)\sqrt{-w^{2}+\cos(\theta)^{2}}}%
\frac{\sin(\theta)^{3}}{\left(  1-w^{2}\right)  ^{3/2}}dw\,d\theta.
\end{align*}
This can be reduced to a single integral:
\begin{align*}
\frac\pi4%
{\displaystyle\int\limits_{0}^{\pi}}
\,\theta\,\tan(\theta)\left[  \cos(\theta)+\sin(\theta)-1\right]  d\theta &
=\frac\pi4\left[  2+\left(  1+\ln(2)\right)  \pi-4G\right] \\
\  &  =2.8708787614233542583742340...
\end{align*}
using the fact that
\begin{align*}
&  \ \
{\displaystyle\int}
\arccos(w)\frac w{\left(  1-w^{2}\right)  ^{3/2}\sqrt{-w^{2}+\cos(\theta)^{2}%
}}dw\\
\  &  =-\frac1{\sin(\theta)^{2}}\left(  \arccos(w)\sqrt{\frac{-w^{2}%
+\cos(\theta)^{2}}{1-w^{2}}}+%
{\displaystyle\int}
\frac{\sqrt{-w^{2}+\cos(\theta)^{2}}}{1-w^{2}}dw\right)
\end{align*}
and
\[
\left.
\begin{array}
[c]{ccc}%
\left.  \arccos(w)\sqrt{\dfrac{-w^{2}+\cos(\theta)^{2}}{1-w^{2}}}\right|
_{w=0}^{\cos(\theta)} &  & \text{if }0\leq\theta\leq\pi/2,\\
\left.  \arccos(w)\sqrt{\dfrac{-w^{2}+\cos(\theta)^{2}}{1-w^{2}}}\right|
_{\cos(\theta)}^{w=0} &  & \text{if }\pi/2\leq\theta\leq\pi
\end{array}
\right\}  =-\frac\pi2\cos(\theta),
\]
\[
\left.
\begin{array}
[c]{ccc}%
{\displaystyle\int\limits_{0}^{\cos(\theta)}}
\dfrac{\sqrt{-w^{2}+\cos(\theta)^{2}}}{1-w^{2}}dw &  & \text{if }0\leq
\theta\leq\pi/2,\\%
{\displaystyle\int\limits_{\cos(\theta)}^{0}}
\dfrac{\sqrt{-w^{2}+\cos(\theta)^{2}}}{1-w^{2}}dw &  & \text{if }\pi
/2\leq\theta\leq\pi
\end{array}
\right\}  =\frac\pi2\left(  1-\sin(\theta)\right)  .
\]
We have not attempted to extend these formulas for $\beta\neq\pi/2$. It is
intriguing that quadrantal triangles ($b=\pi/2$) should present an unevaluated
integral $3.05...\,$while right-angled triangles ($\beta=\pi/2$) give an
integral $2.87...$ expressible in closed-form.

\section{Unconstrained Scenario}

Miles \cite{Ms} proved that
\[
\operatorname*{E}((\alpha+\beta+\gamma-\pi)(a+b+c))=\frac32\pi^{2}-6
\]
where $\alpha+\beta+\gamma-\pi$ is the area $V$ of the spherical triangle and
$a+b+c$ is perimeter $S$. (The notation $V$, $S$ appears to be traditional.)
By preceding correlation results,
\[
3\operatorname*{E}(\alpha\,a)+6\left(  \frac{\pi^{2}}4\right)  -3\pi\left(
\frac\pi2\right)  =3\operatorname*{E}(\alpha\,a)+6\operatorname*{E}%
(\alpha\,b)-3\pi\operatorname*{E}(a)=\frac32\pi^{2}-6
\]
hence $\operatorname*{E}(\alpha\,a)=\pi^{2}/2-2$. It remains to verify Miles' argument.

Up to now, our random spherical triangles have been built using independent
uniform vertices. From now on, they will be built using independent uniform
great circles. By duality, $\operatorname*{E}(V\,S)=3\pi^{2}/2-6$ under either convention.

Let $k$ independent uniform great circles be placed on the unit sphere. The
number of polygonal cells determined is $k^{2}-k+2$ almost always. Randomly
select one of the cells (endowed with equal weighting) and denote the density
for $(V,S)$ by $f_{k}(v,s)$. For example, if $k=2$, then \cite{Gw}
\[%
\begin{array}
[c]{ccc}%
f_{2}(v,s)=\dfrac14\sin\left(  \dfrac v2\right)  \delta(s-2\pi) &  & \text{if
}0\leq v\leq2\pi
\end{array}
\]
and $\delta$ is the Dirac delta function. No formulas for $f_{k}(v,s)$ are
known for $k\geq3$, although when $k=3$ marginal densities for $V$ and for $S$
are well-understood \cite{FJ}.

Let the cells be labeled randomly by the integers $1$, $2$, $3$, $\ldots$,
$k^{2}-k+2$. It is not allowed, for example, to specify that cell $1$ cover
the north pole and that cells $2$, $3$ be adjacent to it. The labeling must be
independent of all features of the tessellation. Hence, for the preceding
experiment, a cell was selected merely by generating a uniform integer
$j\in[1,$ $k^{2}-k+2]$. This is the most basic sampling technique.

We wish to examine alternative methods for selecting a cell. Suppose that the
weighting is proportional to cellular area. Let $C_{j}$ denote the event that
a uniform point falls in cell $j$, where $1\leq j\leq k^{2}-k+2$. If the
volume $V_{j}$ of the cell is $v$, then the probability of $C_{j}$ is
$v/(4\pi)$; unconditionally it is $\operatorname*{E}_{k}(V)/(4\pi)$. The
density for $(V,S)$ here is
\begin{align*}
g_{k}(v,s)  &  =\operatorname*{P}\nolimits_{k}\left\{  \left.  V_{j}%
\in[v,v+dv]\text{ and }S_{j}\in[s,s+ds]\right|  \,C_{j}\right\} \\
&  =\frac{\operatorname*{P}\nolimits_{k}\left\{  V_{j}\in[v,v+dv]\text{ and
}S_{j}\in[s,s+ds]\text{ and}\,C_{j}\right\}  }{\operatorname*{P}%
\nolimits_{k}\left\{  C_{j}\right\}  }\\
&  =\frac{\operatorname*{P}\nolimits_{k}\left\{  C_{j}\left|  V_{j}%
\in[v,v+dv]\text{ and }S_{j}\in[s,s+ds]\right.  \right\}  \,f_{k}%
(v,s)}{\operatorname*{P}\nolimits_{k}\left\{  C_{j}\right\}  }\\
&  =\frac{\left(  v/(4\pi)\right)  f_{k}(v,s)}{\operatorname*{E}_{k}%
(V)/(4\pi)}=\frac{v\,f_{k}(v,s)}{\operatorname*{E}_{k}(V)};
\end{align*}
thus
\begin{equation}
v\,g_{k}(v,s)=\frac{v^{2}\,f_{k}(v,s)}{\operatorname*{E}_{k}(V)}. \label{eq31}%
\end{equation}

Suppose instead that the weighting is proportional to cellular perimeter.
A\ uniform great circle hits $2k$ cells almost always; we then choose one of
these cells at random. Let $C_{j}^{\prime}$ denote the event that a uniform
great circle hits cell $j$ and cell $j$ is subsequently chosen. If the
perimeter $S_{j}$ of the cell is $s$, then the probability of $C_{j}^{\prime}$
is $(s/(2\pi))(1/(2k))$; unconditionally it is $\operatorname*{E}_{k}(S)/(4\pi
k)$. The density for $(V,S)$ here is
\begin{align*}
h_{k}(v,s)  &  =\operatorname*{P}\nolimits_{k}\left\{  \left.  V_{j}%
\in[v,v+dv]\text{ and }S_{j}\in[s,s+ds]\right|  \,C_{j}^{\prime}\right\} \\
&  =\frac{\operatorname*{P}\nolimits_{k}\left\{  V_{j}\in[v,v+dv]\text{ and
}S_{j}\in[s,s+ds]\text{ and}\,C_{j}^{\prime}\right\}  }{\operatorname*{P}%
\nolimits_{k}\left\{  C_{j}^{\prime}\right\}  }\\
&  =\frac{\operatorname*{P}\nolimits_{k}\left\{  C_{j}^{\prime}\left|
V_{j}\in[v,v+dv]\text{ and }S_{j}\in[s,s+ds]\right.  \right\}  \,f_{k}%
(v,s)}{\operatorname*{P}\nolimits_{k}\left\{  C_{j}^{\prime}\right\}  }\\
\  &  =\frac{\left(  s/(4\pi k)\right)  f_{k}(v,s)}{\operatorname*{E}%
_{k}(S)/(4\pi k)}=\frac{s\,f_{k}(v,s)}{\operatorname*{E}_{k}(S)};
\end{align*}
thus
\begin{equation}
v\,h_{k}(v,s)=\frac{v\,s\,f_{k}(v,s)}{\operatorname*{E}_{k}(S)}. \label{eq32}%
\end{equation}

Here is an equivalent definition of $C_{j}$ which is more compatible with that
of $C_{j}^{\prime}$. The intersection of two\ independent uniform great
circles (two diametrically-opposed points $z$ and $-z$) hits two cells almost
always; we then choose one of these cells at random. The new vertex $\pm z$
has four new adjacent cells; upon integrating both sides of (\ref{eq31}), it
becomes clear that
\[
4\operatorname*{E}\nolimits_{k+2}(V)=\frac{\operatorname*{E}_{k}(V^{2}%
)}{\operatorname*{E}_{k}(V)}.
\]
In the same way, with regard to $C_{j}^{\prime}$, the new arc forms the
boundary between two new adjacent cells; upon integrating both sides of
(\ref{eq32}), it becomes clear that
\[
2\operatorname*{E}\nolimits_{k+1}(V)=\frac{\operatorname*{E}_{k}%
(V\,S)}{\operatorname*{E}_{k}(S)}.
\]
Therefore
\[
\frac{\operatorname*{E}_{k-1}(V^{2})}{\operatorname*{E}_{k-1}(V)}%
=2\frac{\operatorname*{E}_{k}(V\,S)}{\operatorname*{E}_{k}(S)}
\]
and, setting $k=3$,
\[
\operatorname*{E}\nolimits_{3}(V\,S)=\frac12\frac{\operatorname*{E}%
\nolimits_{2}(V^{2})}{\operatorname*{E}_{2}(V)}\operatorname*{E}%
\nolimits_{3}(S)=\frac12\frac{2\left(  \pi^{2}-4\right)  }\pi\frac{3\pi
}2=\frac32\pi^{2}-6
\]
as was to be shown.

For $k=3$, the number $N$ of cellular vertices is $3$ almost always. For $k=4
$, the number $N$ is $3$ with probability $4/7$ and $4$ with probability $3/7
$. Recursive equations in $k$ for second order moments of $V$, $S$, $N$ appear
in \cite{Ms, CM} which vastly generalize our discussion here.

\section{Acknowledgement}

I am grateful to Richard Cowan for providing the clearer version of Miles'
proof that appears here. Much more relevant material can be found at \cite{F3,
F4}, including experimental computer runs that aided theoretical discussion here.

\section{Addendum I}

M. Larry Glasser reduced the integral $3.05...\,$to an expression
\[
\frac{\pi^{2}}{2}-\frac{4G}{\pi}-\frac{2}{\pi}\,_{4}F_{3}\left(  \dfrac{1}%
{2},\dfrac{1}{2},1,1;\dfrac{3}{2},\dfrac{3}{2},\dfrac{3}{2};1\right)
\]
where
\[
_{4}F_{3}\left(  \dfrac{1}{2},\dfrac{1}{2},1,1;\dfrac{3}{2},\dfrac{3}%
{2},\dfrac{3}{2};x\right)  =\frac{\sqrt{\pi}}{8}%
{\displaystyle\sum\limits_{n=0}^{\infty}}
\frac{\Gamma(n+1/2)^{2}\Gamma(n+1)^{2}}{\Gamma(n+3/2)^{3}}\frac{x^{n}}{n!}.
\]
He and Jonathan Borwein independently found that
\[
_{4}F_{3}\left(  \dfrac{1}{2},\dfrac{1}{2},1,1;\dfrac{3}{2},\dfrac{3}%
{2},\dfrac{3}{2};1\right)  =%
{\displaystyle\int\limits_{0}^{\pi/2}}
\frac{\operatorname*{Li}\nolimits_{2}(\sin(\theta))-\operatorname*{Li}%
\nolimits_{2}(-\sin(\theta))}{2}d\theta
\]
where $\operatorname*{Li}\nolimits_{2}$ is the dilogarithm function. Let
$\operatorname*{agm}(x,y)$ denote the common limit of sequences $\{a_{n}\}$
and $\{b_{n}\}$ defined via \cite{BB}
\[%
\begin{array}
[c]{cccccccc}%
a_{0}=x, &  & b_{0}=y, &  & a_{n}=\dfrac{a_{n-1}+b_{n-1}}{2}, &  & b_{n}%
=\sqrt{a_{n-1}b_{n-1}} & \text{for }n\geq1.
\end{array}
\]
David Broadhurst's preferred integral for $3.05...\,$is
\[
-%
{\displaystyle\int\limits_{\pi/2}^{\pi}}
\frac{\sin(\theta)+\theta\cos(\theta)}{\operatorname*{agm}(1,\sin(\theta
))}d\theta
\]
because it permits quick high-precision numerical computation.

\section{Addendum II}

Jacopo D'Aurizio \cite{A1} initiated a flurry of research activity, starting
with \cite{A2} and culminating with Vladimir Reshetnikov's striking formula
\cite{A3}:%
\[
_{4}F_{3}\left(  \dfrac{1}{2},\dfrac{1}{2},1,1;\dfrac{3}{2},\dfrac{3}%
{2},\dfrac{3}{2};1\right)  =\frac{3\pi^{3}}{16}+\frac{\pi}{4}\ln
(2)^{2}-4\operatorname{Im}\left(  \operatorname{Li}_{3}(1+i)\right)  .
\]
While the real part of the complex tetralogarithm can be rewritten as
\cite{A4}%
\[
\operatorname{Re}\left(  \operatorname{Li}_{3}(1+i)\right)  =\frac{\pi^{2}%
}{32}\ln(2)+\frac{35}{64}\zeta(3)
\]
the imaginary part has defied all efforts at simplification thus far.

\end{document}